\documentclass[11pt]{article}
\usepackage{graphicx}
\usepackage{amsmath,amsthm,amssymb,amsfonts}
\newtheorem{theo}{Th\'eor`eme}
\newtheorem{prop}{Proposition}
\newtheorem{cor}{Corollaire}
\newtheorem{lem}{Lemme}

\theoremstyle{definition}

\newtheorem{ex}{Exemple}
\newtheorem{remarque}{Remarque}

\newcommand{\RR}{\mathbb{R}}
\newcommand{\ZZ}{\mathbb{Z}}

\title{Sur la g\'eom\'etrie systolique des vari\'et\'es de Bieberbach}

\author{Chady El Mir  \and Jacques Lafontaine}

\begin{document}

\maketitle

\begin{center}

{              Institut de Math\'ematiques et Mod\'elisation de Montpellier \\
CNRS, UMR 5149\\
Universit\'e Montpellier 2\\
              CC 0051, Place Eug\`ene Bataillon\\
              F-34095 Montpellier Cedex 5, France\\
          }
\end{center}

\begin{abstract}
The \emph{systole} of a compact non simply connected Riemannian manifold
is the smallest length of a non-contractible closed curve ; the \emph{systolic ratio}
is the quotient $(\mathrm{systole})^n/\mathrm{volume}$.
Its supremum, on the set of all the riemannian metrics, is known to be finite
for a large class of manifolds, including the $K(\pi,1)$.

We study the optimal systolic ratio of compact, $3$-dimensional non orientable Bieberbach manifolds,
and prove that it cannot be realized by a flat metric.


\vspace{0.5cm}

\end{abstract}

\section{Introduction}

\subsection{Pr\'esentation du r\'esultat}

Un invariant naturel d'une vari\'et\'e riemannienne compacte (non simplement connexe) $(M^n,g)$
est sa  \emph{systole}, not\'ee $\mathrm{sys}(g)$. C'est la plus petite longueur d'une courbe ferm\'ee non contractile. 
Pour avoir un invariant homog\`ene, on introduit le \emph{quotient systolique}
$\frac{\mathrm{sys}(g)^n}{\mathrm{vol}(g)}$.

 Un r\'esultat fondamental de M. Gromov (cf. \cite{gromov}), faisant  suite \`a une propagande
inlassable de M. Berger (voir  \cite{berger}),
assure que \emph{si $M^n$ est \emph{essentielle} la borne sup\'erieure des quotients systoliques,
quand $g$ parcourt l'ensemble des m\'etriques riemanniennes sur $M^n$
est finie.}
Une vari\'et\'e compacte $M$ est dite \emph{essentielle} 
s'il existe une application continue de $M$ dans un $K(\pi,1)$ qui envoie
la classe fondamentale sur une classe non triviale.
Concernant la g\'eom\'etrie systolique, deux r\'ef\'erences incontournables sont
le rapport \cite{bki} de M. Berger et la monographie \cite{katz} de M. Katz.

Les surfaces compactes autres que $S^2$ sont essentielles, et le th\'eor\`eme de Gromov
est une g\'en\'eralisation profonde de r\'esultats ant\'erieurs concernant
le tore $T^2$ (C.~Loewner, voir \cite{ghl} p.95-96 ou \cite{gromov1} p.295-296 pour une preuve),
et le plan projectif  (P.M.~Pu, voir \cite{pu}).

Pour ces deux vari\'et\'es,
le r\'esultat est beaucoup plus pr\'ecis : on a en fait 
$$
\frac{(\mathrm{sys}(g))^2}{\mathrm{aire}(g)}\le \frac{(\mathrm{sys}(g_0))^2}{\mathrm{aire}(g_0)}
$$
o\`u $g_0$ est plate hexagonale (cas du tore)
ou \`a courbure constante positive (cas du plan projectif).
De plus, dans les deux cas, l'\'egalit\'e caract\'erise les m\'etriques 
homoth\'etiques \`a $g_0$.

Il existe un troisi\`eme cas, \'elucid\'e par C.~Bavard dans  \cite{bavard} puis par T.~Sakai dans \cite{sakai} o\`u la borne sup\'erieure
du quotient systolique est connue, et r\'ealis\'ee, celui de la 
bouteille de Klein. Mais
la m\'etrique optimale 
est singuli\`ere,  plus pr\'ecis\'ement $C^1$ par morceaux. Voir \ref{kleinbavard} ci-dessous. Cet exemple
joue ici un r\^ole d\'ecisif.

Pour les vari\'et\'es essentielles de dimension sup\'erieure, on ne connait que tr\`es peu
de choses.  On ignore par exemple, pour des
vari\'et\'es apparemment aussi simples que les tores et les projectifs r\'eels, si les
m\'etriques \`a courbure constante sont optimales.
L'argument de Gromov, puissant mais tr\`es g\'en\'eral,
ne donne qu'une majoration tr\`es grossi\`ere du quotient systolique. A fortiori, on a tr\`es peu d'informations
sur les m\'etriques qui pourraient r\'ealiser sa borne sup\'erieure.
Ivan Babenko nous a communiqu\'e le r\'esultat suivant :

\emph{si $g$ est une m\'etrique optimale (\'eventuellement singuli\`ere) sur $M$, les 
courbes systoliques (c'est-\`a-dire les courbes ferm\'ees non contractiles de longueur $\mathrm{sys}(g)$) recouvrent $M$.}

Cette propri\'et\'e est satisfaite par les tores plats, et par les projectifs r\'eels munis de leur m\'etrique \`a courbure constante.
Sur les vari\'et\'es compactes de dimension $3$ qui portent une m\'etrique plate, les m\'etriques qui
optimisent le quotient systolique \emph{parmi les m\'etriques plates}
v\'erifient  aussi cette propri\'et\'e, \`a l'exception d'une d'entre elles.
(Il s'agit de la vari\'et\'e dont le  
groupe fondamental est le groupe appel\'e $G_6$ par J.A. Wolf, cf. \cite{wolf} p.117-118, et le dessin
tr\`es  suggessif de \cite{th}, p.236 ; ce groupe est engendr\'e par trois sym\'etries gliss\'ees par rapport \`a trois droites
deux \`a deux
orthogonales).

Dans ce travail, nous nous int\'eressons aux \emph{vari\'et\'es de Bieberbach}, c'est--\`a--dire
aux vari\'et\'es compactes qui portent une m\'etrique riemannienne plate.
Ces vari\'et\'es sont des $K(\pi,1)$, donc le th\'eor\`eme de Gromov s'applique.
L'exemple de la bouteille de Klein, et celui de la vari\'et\'e de dimension $3$ que nous venons d'\'evoquer,
conduisent \`a penser que,
au moins  
 pour celles de ces vari\'et\'es qui ne sont pas des tores,
il ne faut pas s'attendre \`a ce que les m\'etriques plates soient optimales. 
Notre r\'esultat est le suivant.

\emph{%
Soit $M$ une vari\'et\'e de Bieberbach de dimension $3$ non orientable.
Alors il existe sur $M$ une m\'etrique riemannienne $g$ telle que,
pour toute m\'etrique plate $h$,
$$
\frac{(\mathrm{sys}(h))^3}{\mathrm{aire}(h)} <   \frac{(\mathrm{sys}(g))^3}{\mathrm{vol}(g)}
$$
}

\subsection{R\'esum\'e de la preuve}
 Nous commen\c cons par v\'erifier (cf. \ref{plat}) 
que nos m\'etriques plates s'obtiennent par des suspensions de la bouteille de Klein.
Dans \ref{sing}, nous \'etudions d'un peu plus pr\`es la m\'etrique singuli\`ere
sur la bouteille de Klein d\'ecouverte par C.~Bavard,
et en particu\-lier ses isom\'etries. Celles-ci permettent de construire par
suspension des m\'etriques singuli\`eres en dimension $3$. Le calcul direct des systoles,
fait dans \ref{fin}, 
met en \'evidence des m\'etriques singuli\`eres pour lesquelles
le quotient systolique
est sup\'erieur \`a celui de toute m\'etrique plate sur la m\^eme vari\'et\'e.

Pour conclure, il suffit d'approcher la m\'etrique singuli\`ere obtenue par une m\'e\-trique lisse, et d'utiliser
la continuit\'e de la systole pour la topologie $C^0$ (voir les d\'etails dans  \ref{continuite}).

\subsection{Notations}
Si $V$ est un sous-espace affine d'un espace affine euclidien $E$, 
et $a$ un vecteur parall\`ele  \`a $V$, 
on note $\sigma_{V,a}$ la sym\'etrie gliss\'ee
produit de la r\'eflexion orthogonale par rapport \`a $V$
et de la translation de vecteur $a$.

Quand une isom\'etrie de $\mathbb{R}^n$ passe au quotient,
nous la noterons par le m\^eme symbole.

On appelle \emph{d\'eplacement} d'une isom\'etrie $\gamma$ d'un espace m\'etrique (qui sera ici une vari\'et\'e
riemannienne, \'eventuellement  singuli\`ere), le nombre 
$$\mathbf{d}(\gamma)=\inf_{p\in M}\mathrm{dist}\big(p,\gamma(p)\big).$$
Rappelons que si $(M,g)$ est une vari\'et\'e riemannienne compacte,
de rev\^etement riemannien $(\tilde M,\tilde g)$,
le groupe fondamental op\`ere sur  $(\tilde M,\tilde g)$ par isom\'etries,
et 
$$\mathrm{sys}\big(M,g\big)=\inf_{\gamma\in \pi_1(M)}\mathbf{d}(\gamma)$$
Ce r\'esultat s'\'etend imm\'ediatement aux m\'etriques singuli\`eres consid\'er\'ees ci-dessous.

\section{Rappels et compl\'ements sur les vari\'et\'es plates}
\label{plat}
\subsection{Vari\'et\'es plates}


Les vari\'et\'es compactes et plates sont les quotients $\mathbb{R}^n/\Gamma$ o\`u $\Gamma$ est un sous groupe  discret $\Gamma$, cocompact et sans points fixes d'isom\'etries affines de $\mathbb{R}^n$. 
Le th\'eor\`eme de Bieberbach
assure que $\Gamma$ est une extension  d'un groupe fini $G$  par un r\'eseau $\Lambda$ de $\mathbb{R}^n$.
Ce r\'eseau est le sous-groupe des \'el\'ements de $\Gamma$ qui sont des translations.
Nous l'appelerons dans la suite le \emph{r\'eseau associ\'e} \`a $\Gamma$.
On a une suite exacte
$$
0\longrightarrow\Lambda\longrightarrow\Gamma\longrightarrow G\longrightarrow e
$$
(En fait, ce r\'esultat se d\'emontre sans supposer que $\Gamma$ op\`ere librement).
Dans la situation qui nous int\'eresse, la vari\'et\'e est le quotient
du tore plat $\mathbb{R}^n/\Lambda$ par un groupe d'isom\'etries isomorphe \`a $G$.

Deux vari\'et\'es compactes plates $\RR^n/\Gamma$ et $\RR^n/\Gamma^\prime$
sont hom\'eomorphes si et seulement si les groupes $\Gamma$ et $\Gamma^\prime$
sont isomorphes. Ces groupes sont alors conjugu\'es par un isomorphisme affine de $\RR^n$, 
ce qui montre  que deux vari\'et\'es compactes plates hom\'eomorphes sont affinement diff\'eomorphes.
Ces vari\'et\'es seront isom\'etriques si et seulement si $\Gamma$ et $\Gamma^\prime$ sont
conjugu\'es par une isom\'etrie affine.

Un proc\'ed\'e classique de construction de vari\'et\'es est la \emph{suspension}
par un diff\'eo\-mor\-phisme $\varphi$ d'une vari\'et\'e $K$ : il
s'agit du quotient de $K\times\mathbb{R}$ par le groupe engendr\'e par
$(p,t)\mapsto \big(\varphi(p), t+1\big)$.
Si maintenant $(K^n,g)$ est une vari\'et\'e riemannienne de dimension $n$ et si $\varphi\in\mathrm{Isom}(K^n,g)$,
le quotient \emph{riemannien} de   $K\times\mathbb{R}$ par 
le groupe engendr\'e par
$(p,t)\mapsto \big(\varphi(p), t+a\big)$, o\`u $a$ est un param\`etre $>0$,
d\'efinit sur cette suspension une m\'etrique riemannienne (bien s\^ur plate si $g$ est plate).
Cette construction sera utilis\'ee syst\'ematiquement dans \ref{plat3} puis dans \ref{fin}.

\subsection{Folklore kleinien} \label{klein}

Le plan euclidien \'etant rapport\'e \`a une base orthonorm\'ee,
les bouteilles de Klein  plates  sont les vari\'et\'es $\mathbf{R}^2/\Gamma$,
o\`u $\Gamma$ est engendr\'e par la sym\'etrie gliss\'ee
$(x,y)\mapsto (x+\frac{a}{2},-y)$
et la translation $(x,y)\mapsto (x,y+b)$.
Nous noterons $K_{a,b}$ (ou simplement $K$ quand la donn\'ee de $a$ et $b$ est sous-entendue ou
inutile), la vari\'et\'e riemanniennne plate ainsi obtenue.

Rappelons que les g\'eod\'esiques ``horizontales" de $K_{a,b}$ sont ferm\'ees de longueur $a$,
\`a l'exception de deux d'entre elles qui sont de longueur $\frac{a}{2}$ (cf. \cite{ghl} p.82-83). Cela permet de voir
$K_{a,b}$ comme le recollement de deux rubans de M\"obius plats de largeur $\frac{b}{2}$ le long de leur bord.
Cette identification n'est pas seulement topologique mais aussi riemannienne : la r\'eflexion 
orthogonale par rapport au bord commun est une isom\'etrie qui \'echange les deux rubans.

La composante neutre de $\mathrm{Isom}(K)$ est form\'ee des translations horizontales
$r_\alpha :(x,y)\mapsto (x+\alpha,y)$ ($\alpha$ \'etant pris modulo $a$).
A signaler : la translation $(x,y)\mapsto (x+\frac{a}{2},y)$
qui co\"\i ncide 
avec
$(x,y)\mapsto (x,-y)$.
Elle laisse fixes les deux g\'eod\'esiques hori\-zontales courtes $y=0$ et $y=b/2$.
C'est la sym\'etrie orthogonale par rapport \`a ces \emph{deux} g\'eod\'esiques simultan\'ement.

Le quotient $\mathrm{Isom}(K)/\mathrm{Isom}_0(K)$ est isomorphe au
groupe du matelas, appel\'e aussi groupe de Klein (c'est plus qu'une co\"\i ncidence !).

Les trois \'el\'ements non triviaux de ce quotient peuvent \^etre repr\'esent\'es par

\begin{enumerate}
\item une r\'eflexion par rapport \`a une g\'eod\'esique verticale, qui
est aussi une sym\'etrie par rapport \`a un point d'une g\'eod\'esique horizontale courte
(en fait, une telle transformation, que nous noterons $S_1$, est \`a la fois une r\'eflexion par
rapport \`a deux g\'eod\'esiques verticales et une sym\'etrie par rapport
\`a deux points distincts).

\item une sym\'etrie par rapport \`a un point du bord commun des deux rubans.
Nous noterons $S_2$ une telle sym\'etrie.

\item la r\'eflexion par rapport \`a ce bord commun, ou la transformation
obtenue par passage au quotient de 
$(x,y)\mapsto (x,y+ b/2)$. Nous noterons $T$ cette derni\`ere transformation.
\end{enumerate}

Notons enfin qu'un diff\'eomorphisme affine d'une bou\-teille de Klein plate est une isom\'etrie.

\subsection{Une description g\'eom\'etrique des vari\'et\'es plates non orientables de dimension $3$}
\label{plat3}

Nous allons voir que toute vari\'et\'e compacte plate non orientable de dimension $3$
peut s'obtenir par suspension d'une bouteille de Klein plate par une isom\'etrie.

Cela se voit en interpr\'etant convenablement (et plus g\'eom\'etriquement)
la description faite par J. Wolf (voir \cite{wolf}, p.120-123, dont nous suivons
les notations) des vari\'et\'es de Bieberbach de dimension $3$. Voir aussi \cite{th}, 4.2. 

i) Type $B_1$. On se donne un plan $P$, deux vecteurs ind\'ependants 
$a_1$ et $a_2$
de $P$, et un vecteur $a_3$ orthogonal \`a $P$.  Le groupe $\Gamma$ est engendr\'e par la sym\'etrie gliss\'ee $\sigma_{P,a_1/2}$
et les translations
de vecteurs $a_2$ et $a_3$.
Notons que $\Gamma$  est d'indice $2$ dans le r\'eseau $\Lambda$ engendr\'e par
$a_1$, $a_2$ et $a_3$, et que 
la vari\'et\'e est diff\'eomorphe au produit d'une bouteille de Klein et d'un cercle.
Les diff\'erents types m\'etriques s'obtiennent en faisant varier les vecteurs
$a_i$ avec les contraintes ci-dessous.

Pour r\'ealiser ces m\'etriques commme des suspensions, on \'ecrit 
 $a_2=\lambda a_1+v$, o\`u $v$ est parall\`ele \`a $P$ et orthogonal \`a $a_1$, et on 
v\'erifie que $\RR^3/\Gamma$ 
 est le quotient riemannien
de
$K_{\vert a_1\vert ,\vert a_3\vert     }\times\mathbb{R}$
par le groupe d'isom\'etries engendr\'e par
$$(p,t)\mapsto \big(r_\alpha(p), t+\vert v\vert\big),\quad \hbox{avec}\
\alpha ={\lambda\vert a_1\vert}.$$


ii) Type $B_2$. On se donne deux  plans parall\`eles $P_1$ et $P_2$, et deux vecteurs
in\-d\'e\-pen\-dants $a_1$ et $a_2$ parall\`eles \`a ces  plans. Le groupe
$\Gamma$ est engendr\'e par
les sym\'etries gliss\'ees $\sigma_{P_1,a_1/2}$ et  $\sigma_{P_2,a_2/2}$.
Soit $\mathbf{k}$ un vecteur unitaire orthogonal \`a ces plans,
et $d\mathbf{k}$ le vecteur de la translation qui envoie $P_1$ sur $P_2$.
Alors  $\sigma_{P_2,a_2/2}\circ \sigma_{P_1,a_1/2}$
est la translation de vecteur $a_3= \frac{a_1+a_2}{2}+ 2d\mathbf{k}$. Le r\'eseau $\Lambda$ associ\'e \`a $\Gamma$
est engendr\'e par $a_1$, $a_2$ et $a_3$.
Les diff\'erents types m\'etriques s'obtiennent en faisant varier $a_1$, $a_2$ et $d=\mathrm{dist}(P_1,P_2)$.

Cette vari\'et\'e est diff\'eomorphe \`a la suspension d'une bouteille de Klein par $T$.
En effet, le sous-groupe  engendr\'e par $\sigma_{P_1,a_1/2}$
et la translation de vecteur $4d\mathbf{k}$ est distingu\'e  dans $\Gamma$ ; il est isomorphe au groupe 
fondamental de la bouteille de Klein, et le quotient est isomorphe \`a $\ZZ$.

Pour expliciter cette identification au niveau des m\'etriques, on \'ecrit comme plus haut
$a_2=\lambda a_1+v$, avec  $v$ orthogonal \`a $a_1$, et on v\'erifie
que $\RR^3/\Gamma$ 
 est le quotient riemannien
de
$K_{\vert a_1\vert ,4d}\times\RR$ par le groupe d'isom\'etries engendr\'e par
$$(p,t)\mapsto \Big(r_{\alpha}\big(T(p)\big), t+\frac{\vert v\vert}{2}\Big),\quad
\hbox{avec}\ \alpha= \vert a_1\vert\frac{\lambda -1}{2}.$$

iii) Type $B_3$. On se donne un plan $P$, une droite $D$ de $P$, et une
 base orthogonale $(a_1,a_2, a_3)$, telle que $a_1$ et  $a_2$ soient parall\`eles \`a $ D$ et  
 $P$ respectivement.
 Le groupe
 $\Gamma$ est engendr\'e par $\sigma_{D,a_1/2}$, $\sigma_{P,a_2/2}$
 et la translation de vecteur $a_3$.
On peut aussi obtenir cette m\'etrique   par suspension de la bouteille de 
Klein :
c'est le quotient riemannien de $K_{\vert a_2\vert  ,\vert a_3\vert }\times \mathbb{R}   $
par le groupe d'isom\'etries engendr\'e par
$$(p,t))\mapsto \big(S_1(p), t+ \vert a_1\vert/2\big)\quad \hbox{(voir \ref{klein})}.$$

iv) Type $B_4$. On se donne un plan $P$, une droite $D$ parall\`ele \`a
$P$ mais non incluse dans $P$,  deux vecteurs orthogonaux
$a_1$ et $a_2$  parall\`eles \`a $D$ et $P$ respectivement.
Le groupe $\Gamma$ est engendr\'e par $\sigma_{D,a_1/2}$ et $\sigma_{P,a_2/2}$.

Le r\'eseau associ\'e \`a $\Gamma$ est engendr\'e par la base orthogonale
$(a_1,a_2,a_3)$, avec $\vert a_3\vert =4\mathrm{dist}(P,D)$.
Cette fois, la m\'etrique est le quotient riemannien de
 $K_{\vert a_2\vert,\vert a_3\vert}\times \mathbb{R}$
par le groupe d'isom\'etries engendr\'e par
$$(p,t))\mapsto \big(S_2(p), t+ \vert a_1\vert/2\big)\quad\hbox{(voir encore \ref{klein})}.$$

\begin{remarque} Cette discussion montre que 
chaque classe de diff\'eo\-morphisme de va\-ri\'et\'e plate com\-pacte non orientable de dimension $3$
est associ\'ee \`a un \'e\-l\'e\-ment de  du groupe des transformations affines de $K$ quotient\'e
par sa composante neutre.
\end{remarque}

\section{M\'etriques singuli\`eres}
\label{sing}
\subsection{Pr\'esentation de la bouteille de Klein--Bavard}
\label{kleinbavard}
Nous travaillerons avec des m\'etriques riemanniennes \emph{singuli\`eres} dans le sens suivant.

\begin{enumerate}
\item  Elles sont continues, c'est-\`a-dire que les coefficients $g_{ij}$, exprim\'es
en cartes locales, sont continus.
\item La vari\'et\'e est une r\'eunion de domaines \`a bord d'int\'erieurs deux \`a deux disjoints, tels qu'\`a l'int\'erieur de chaque domaine la m\'etrique est lisse.
\end{enumerate}

Une telle m\'etrique d\'efinit un espace de longueur
(voir \cite{bbi} pour des d\'etails sur cette notion).
On peut donc parler de g\'eod\'esiques. De plus, la mesure riemannienne se d\'efinit
exactement comme dans le cas lisse.

Pour les m\'etriques qui nous int\'eressent, le mod\`ele local est le suivant.
On part de la sph\`ere ronde, et on rep\`ere les points par la latitude
$\phi$ et la longitude $\theta$. Pour $\phi_o\in ]0,\pi/2[$, 
soit $\Sigma_{\phi_o}$ le domaine d\'efini par $\vert\phi\vert \le \phi_o$.
Sur  $\Sigma_{\phi_o}$, la m\'etrique ronde est donn\'ee par
$d\phi^2 +\cos^2\phi d\theta^2$. A partir de l\`a,
on introduit  sur
$\RR^2$ la m\'etrique riemannienne singuli\`ere (au sens pr\'ec\'edent) 
\begin{equation}
d\phi^2 +f^2(\phi)d\theta^2,
\label{sph}
\end{equation}
 o\`u $f$ est la fonction 
$2\phi_0$-p\'eriodique qui co\"\i ncide avec $\cos\phi$ sur l'intervalle
$[-\phi_o, \phi_o]$ (comparer \`a \cite{bavard}).

\begin{ex} La m\'etrique sur la bouteille de Klein qui donne le quotient systolique
optimal s'obtient pour $\phi_o=\frac{\pi}{4}$ en quotientant 
le plan muni de la m\'etrique \ref{sph}
par l'action du groupe
engendr\'e par
$$ (\theta,\phi)\mapsto (\theta+\pi,-\phi)\quad\hbox{et}
\quad  (\theta,\phi)\mapsto (\theta,\phi+ 4\phi_0).$$
\end{ex}

Nous noterons $(\mathbf{K},b)$ la bouteille de Klein munie de cette m\'etrique,
$(\mathbf{T}^2, b)$ son rev\^etement orientable.
Il se trouve que $(\mathbf{K},b)$ a le m\^eme groupe d'iso\-m\'e\-tries que celui d'une bouteille de Klein plate
(attention, la propri\'et\'e analogue pour $(\mathbf{T}^2,b)$ est fausse !).
Cela peut se voir par un calcul direct, mais il sera utile de voir les choses
g\'eom\'etriquement :
le domaine sph\'erique \`a bord $\Sigma_{\phi_o}$
vu plus haut peut se voir comme 
$$ \{(x,y,z)\in\RR^3, x^2+y^2+z^2=1, \vert z\vert\le  \sin\phi_o\}.$$
Le cylindre obtenu en prenant tous les tranlat\'es de vecteur $2n\sin \phi_0\vec{k}$
de $\Sigma_{\phi_o}$, muni de la structure de longueur induite par la m\'etrique euclidienne,
est un rev\^etement riemannien de $(\mathbf{T}^2,b)$ (et donc de $(\mathbf{K},b)$).
Nous le noterons $(\mathbf{C},b)$. 

On r\'ecup\`ere $(\mathbf{K},b)$ en quotientant $(\mathbf{C},b)$  par le groupe engendr\'e par les antipodies 
par rapport aux centres des sph\`eres portant ces domaines  (deux domaines contigus suffisent).
On notera $\sigma$ l'une quelconque de ces antipodies.

Alors les isom\'etries suivantes de ce cylindre passent au quotient. Si on adopte la formulation analytique de (\ref{sph}),
 elles sont les exactes analogues des isom\'etries d\'ecrites en \ref{klein}.
\begin{enumerate}
\item les rotations autour de l'axe (analogues aux translations ``horizontales").
\item les r\'eflexions par rapport \`a un plan m\'eridien (analogues aux r\'eflexions par rapport
\`a une g\'eod\'esique ``verticale", c'est-\`a-dire aux transformations de type $S_1$).
\item les r\'eflexions par rapport \`a un diam\`etre d'un cercle singulier
(analogues aux  sym\'etries par rapport \`a un point du bord commun des deux rubans,
c'est-\`a-dire aux transformations de type $S_2$).
\item la r\'eflexion par rapport \`a un plan d'un cercle singulier (qui donne par passage au quotient
une transformation \'equivalente (modulo $\mathrm{Isom}_o$)  \`a celle obtenue en faisant passer au quotient la translation
de vecteur $2\sin \phi_0\vec{k}$. Cette derni\`ere transformation est l'ana\-logue de $(x,y)\mapsto (x,y+\frac{b}{2})$.
\end{enumerate}

\subsection{Quelques propri\'et\'es g\'eom\'etriques de la bouteille singuli\`ere}

\begin{prop} \label{geod} Localement, une g\'eod\'esique de $(\mathbf{C},b)$, $(\mathbf{T}^2,b)$
ou $(\mathbf{K},b)$ est de l'un des types suivants :

1)  un arc de grand cercle (\'eventuellement tangent \`a une ligne singuli\`ere).

2) la juxtaposition d'un arc de grand cercle 
tangent \`a une ligne singuli\`ere, et d'un arc de cette ligne singuli\`ere.

3) la juxtaposition de deux arcs de grand cercle ayant une extr\'emit\'e commune $r$
 sur une ligne singuli\`ere, localement de part et d'autre de cette ligne,
qui fait
avec ces deux arcs des angles dont la somme est $\pi$.

\end{prop}
\begin{proof}[Preuve]
Il suffit de juxtaposer les remarques suivantes :

\begin{enumerate}
\item
 Les courbes contenues dans l'int\'erieur d'un $\Sigma$  sont 
des g\'eod\'esiques si et seulement si ce sont des arcs de grand cercle, d'apr\`es la th\'eorie classique.
\item
Une  courbe contenue dans un $\Sigma$ et dont les extr\'emit\'es $p$ et $q$
sont sur la ligne singuli\`ere est plus longue  que le plus court des deux arcs de la
ligne singuli\`ere qui joignent $p$ et $q$. 
Pour le voir, il suffit d'expliciter la longueur en utilisant la formule (\ref{sph}) et de remarquer
que la fonction $\cos\phi$ est d\'ecroissante.
\item
Soit $c$ une courbe $C^1$ par morceaux, param\'etr\'ee par longueur d'arc,
qui est la juxtaposition d'une courbe $c_1$ contenue dans un $\Sigma$, ayant une extr\'emit\'e $r$
sur une ligne singuli\`ere, et d'une courbe $c_2$ partant de $r$ et contenue dans cette ligne singuli\`ere.
Si $c_1^\prime(r)\not=c_2^\prime(r)$, d'apr\`es la formule de variation premi\`ere,
$c$ peut \^etre raccourcie en gardant les m\^emes extr\'emit\'es.
\item  
soit enfin une courbe dont les extr\'emit\'es $p$ et $q$ sont dans deux exemplaires contigus
de  $\Sigma$, et qui intersecte la ligne singuli\`ere en un point. 
Si elle minimise la longueur des courbes joignant $p$ et $q$, elle est n\'e\-cessai\-rement du troisi\`eme type d\'ecrit dans
l'\'enonc\'e en raison de la formule de variation premi\`ere.
\end{enumerate}
\end{proof}

\begin{remarque}  A partir de l\`a, on voit que deux points appartenant \`a des domaines $\Sigma$ contigus sont
joints en g\'en\'eral par un segment  g\'eod\'esique unique, \'eventuellement par deux.
Ce dernier cas n'arrive que pour des points ayant des latitudes oppos\'ees : la r\'eflexion orthogonale
par rapport au plan d\'efini par l'axe des $z$ et les deux points \'echange alors ces segments.
\end{remarque}

\begin{figure}
\begin{center}
\includegraphics[scale=0.7]{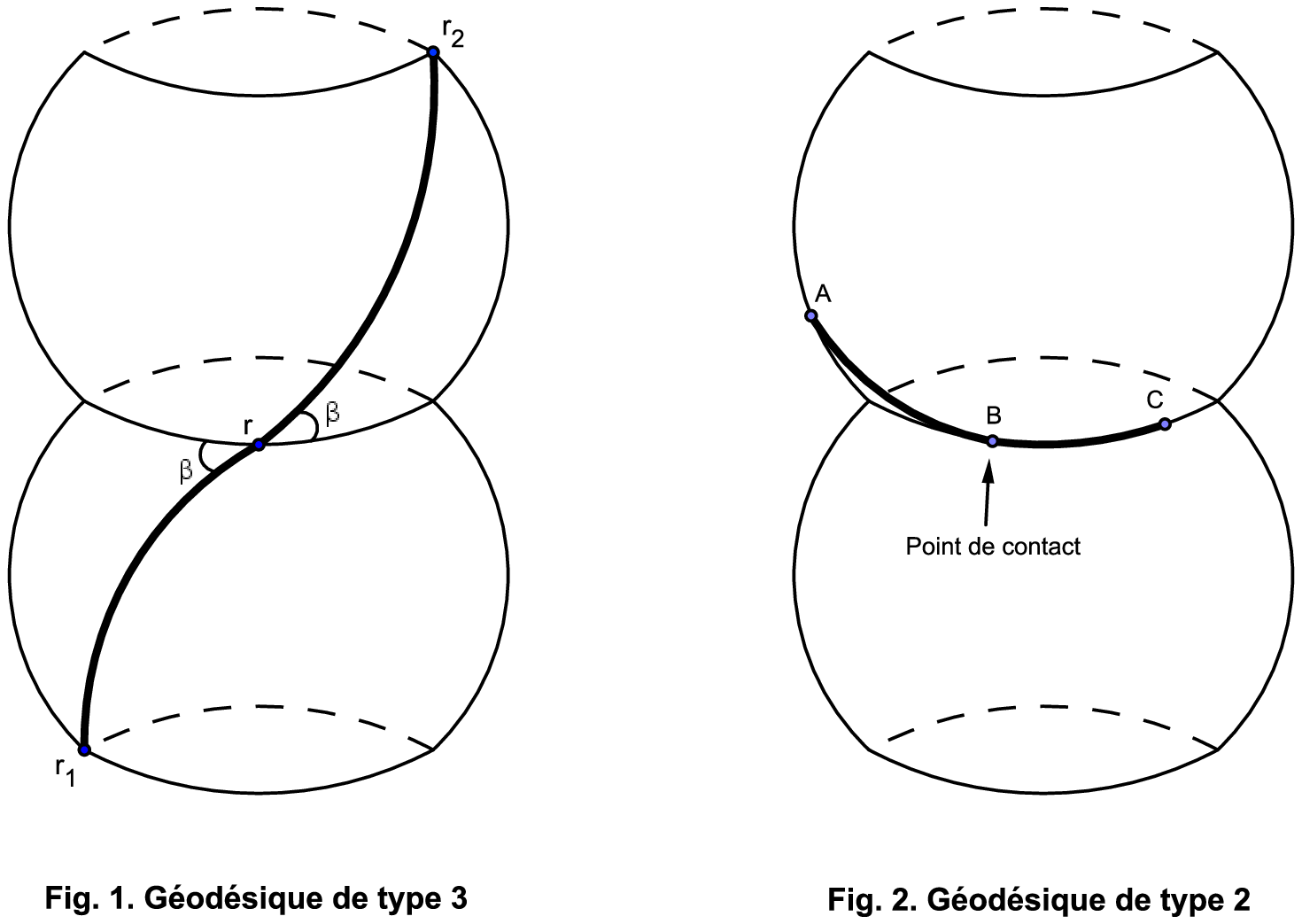}
\end{center}
\end{figure}
 
D\'esignons par $T_{\delta}$
le d\'eplacement h\'elico\"\i dal obtenu en composant la translation
de vecteur $2\sin \phi_0\vec{k}$ et la rotation d'angle $\delta$ autour de l'axe des $z$,
\emph{consid\'er\'e comme isom\'etrie du cylindre $(\mathbf{C},b)$}.
Nous avons vu que cette isom\'etrie passe au quotient pour $(\mathbf{T},b)$
comme pour  $(\mathbf{K},b)$.

\begin{prop}\label{memedist}
Sur  $(\mathbf{T},b)$, la fonction $\mathrm{dist}\big(p,(T_{\delta}(p)\big)$
est constante.
\end{prop}
\begin{proof}[Preuve]
On v\'erifie facilement qu'il suffit de faire la preuve pour  $(\mathbf{C},b)$.
On remarque que $T_{\delta}$ est la restriction au cylindre d'un d\'eplacement de $\RR^3$ 
qui est le produit de deux r\'eflexions
par rapport 
\`a des  droites orthogonales \`a l'axe des $z$,
faisant entre elles un angle \'egal \`a $\delta/2$, et de plus courte distance
(dans l'espace euclidien) $\sin\phi_o$.
Soit $\gamma$ une g\'eod\'esique minimisante joignant $p$  et 
$T_{\delta}(p)$. Admettons un instant que cette g\'eod\'esique rencontre 
une ligne singuli\`ere en un unique point $r$, et que l'on se trouve donc dans la troisi\`eme
situation de la proposition \ref{geod}.
Cet arc est inclus dans la juxtaposition de deux arcs de grand cercle $r_1r$ et $rr_2$,
contenus dans deux domaines sph\'eriques contigus $\Sigma_1$ et $\Sigma_2$, et faisant le m\^eme
angle $\beta$ avec les lignes singuli\`eres en $r$, $r_1$, $r_2$ (voir figure 1). Notons aussi
que $\gamma$ se prolonge en une g\'eod\'esique de $(\mathbf{C},b)$
form\'ee de tels arcs de grand cercle, et que $T_{\delta}(r_1)=r$,
$T_{\delta}(r)=r_2$. Choisissons les r\'eflexions dont la compos\'ee est $T_\delta$ comme suit : 
l'axe de la premi\`ere passant par le milieu de l'arc $r_1r$, l'axe de la seconde passant par $r$.
Cela permet de voir que $T_\delta$ laisse la g\'eod\'esique qui prolonge $\gamma$
globalement invariante. Comme il n'y a pas de point fixe, c'est une translation le long de cette
g\'eod\'esique. L'assertion de l'\'enonc\'e en r\'esulte, puisque toute la situation 
est invariante par rotation autour de l'axe des $z$.

Reste \`a justifier notre affirmation sur $\gamma$. Elle est clairement vraie pour
$\delta=0$ et donc pour $\delta$ assez petit. Dans ce cas, un calcul explicite 
que nous laissons au lecteur montre que
$$\tan\delta=\frac{-2\sqrt{2}\cot\beta}{\cot^2\beta-2},$$
et donc 
que tous les $\delta$ possibles peuvent \^etre obtenus ainsi. On applique alors
la remarque qui suit la proposition pr\'ec\'edente.
\end{proof}

\emph{Dor\'enavant, nous suppposerons que $\phi_0=\frac{\pi}{4}$.}

Pour les rotations $r_\alpha$ on a la propri\'et\'e suivante.

\begin{prop}\label{minsing}
Soit $p$ un point du tore $(\mathbf{T},b)$ 
et $p_0$ un point d'une ligne singuli\`ere. Alors
$$\mathrm{dist}_{(\mathbf{T},b)}(p,r_{\alpha}(p))\ge \mathrm{dist}_{(\mathbf{T},b)}(p_0,r_{\alpha}(p_0))$$
pour tout $\alpha\in[0,\pi]$
\end{prop}
Ainsi, le minimum est obtenu pour les points de la ligne 
singuli\`ere.

\begin{proof}[Preuve]
On utilise la proposition \ref{geod}, en remarquant qu'il suffit de faire le calcul dans 
$(\mathbf{C},b)$.
Soit $\Sigma$ le domaine sph\'erique o\`u se trouve $p$ (et donc $q=r_\alpha(p)$).
Tant que le grand cercle passant par ces deux points est 
enti\`erement dans $\Sigma$ (ce qui est le cas si leur latitude commune $\beta$ est assez petite), la formule
fondamentale de la trigonom\'etrie sph\'erique donne
$$\mathrm{dist}\big(p,r_\alpha(p)\big)=\arccos(\sin^2\beta+\cos^2\beta\sin\alpha),$$  et
cette distance est une fonction d\'ecroissante de $\beta$.

Si ce n'est pas le cas, on est dans le deuxi\`eme cas de la Proposition \ref{geod}.
Le segment g\'eod\'esique qui joint $p$ et $q$
est form\'e de deux  arcs de grand cercle ayant une extr\'emit\'e en $p$ (resp. $q$), et tangents
en leur autre extr\'emit\'e, que nous noterons $p^\prime$ et $q^\prime$,
\`a une ligne singuli\`ere. 
Si $\psi$ est la diff\'erence de longitude entre $p$ et  $p^\prime$ 
(ou entre $q$ et   $q^\prime$) on v\'erifie que $\cos\psi=\tan \beta$.
On a alors
$$\mathrm{dist}(p,p^\prime)= \mathrm{dist}  (q,q^\prime)=\arccos(\frac{\cos\beta\cos\psi}{\sqrt 2}+\frac{\sin\beta}{\sqrt 2})=\arccos(\sqrt 2\sin\beta)$$
et $$\mathrm{dist}(p^\prime,q^\prime)=\frac{\alpha-2\psi}{\sqrt 2}$$
(on v\'erifie au passage que l'on est dans cette situation si $\beta\ge\arctan(\cos\alpha/2$).

Finalement dans ce cas
$$\mathrm{dist}\big((p,r_\alpha(p)\big)=\frac{\alpha}{\sqrt 2}-{\sqrt 2}\arccos(\tan\beta)+2\arccos({\sqrt 2}  \sin\beta).$$
La fonction de $\beta$ ainsi obtenue a pour d\'eriv\'ee
$$\frac{\sqrt 2}{\cos\beta\sqrt{1-2\sin^2\beta}}(1-2\cos^2\beta).$$
Les expressions de $\mathrm{dist}\big((p,r_\alpha(p)\big)$ pour $\beta\in[0,\beta_0]$, puis pour
$\beta\in [\beta_0,\pi/4]$ montrent que l'on a une fonction d\'ecroissante de $\beta$.
Son minimum est atteint quand $p$ est sur une ligne singuli\`ere. Il vaut $\alpha\sqrt 2$.
\end{proof}

En passant au quotient, on d\'eduit de ce qui pr\'ec\`ede les r\'esultats analogues
pour la bouteille de Klein.

\begin{cor} \label{dklein} Pour $(\mathbf{K},b)$, on a
$$\mathbf{d}(r_\alpha))=\inf\big\{\frac{\alpha}{\sqrt 2}, \pi -\alpha\big\}\quad
\mathrm{et}\quad\mathbf{d}(T_\delta)=
\inf\big\{\frac{\pi -\delta}{\sqrt 2}, \arccos\frac{\cos\delta-1}{2}\big\}$$
\end{cor}
\begin{proof}
En fait, 
$$\mathbf{d}(r_\alpha)=\inf_{p\in\Sigma}\{\mathrm{dist}\big(p,r_\alpha(p)\big),\mathrm{dist}\big(\sigma(p),r_\alpha(p)\big)\}.$$
Mais comme les points $\sigma(p)$ et $r_\alpha(p)$ sont situ\'es de part et d'autre de l'\'equateur de $\Sigma$,
le segment g\'eod\'esique qui les joint est contenu dans $\Sigma$ (premier cas de la proposition
\ref{geod}). On v\'erifie alors que $\mathrm{dist} \big(\sigma(p),r_\alpha(p)\big)$, est une fonction
d\'ecroissante de la latitude, d'o\`u le r\'esultat pour $r_\alpha$.

Pour \'etudier  $\mathbf{d}(T_\delta)$, il suffit de consid\'erer deux exemplaires contigus $\Sigma_1$
et $\Sigma_2$, et de calculer l'infimum analogue au pr\'ec\'edent 
quand $p$ parcourt $\Sigma_1$. On sait d\'ej\`a d'apr\`es la Proposition \ref{memedist}
que $\mathrm{dist}\big(p,T_\delta(p)\big)$ est constante (et \'egale \`a 
$\arccos\frac{\cos\delta-1}{2}$). Passons \`a $\mathrm{dist}\big(T_\delta(p), \sigma(p)\big)$.
Soit $c$ le segment g\'eod\'esique qui joint ces points.
On est dans le deuxi\`eme
ou le troisi\`eme cas de la Proposition \ref{geod}.
Dans tous les cas, 
soit
$r$ son premier point de rencontre avec la ligne singuli\`ere, et
 $q$ le sym\'etrique de $T_\delta(p)$ par rapport \`a cette ligne singuli\`ere.

Dans $\Sigma_1$, on passe de $q$ \`a $\sigma(p)$ par une rotation d'angle $\pi-\delta$.
Par ailleurs, ces deux points sont reli\'es par 
la  courbe $\tilde c$ obtenue en prenant la sym\'etrique de la partie de $c$ de $T_\delta(p)$ \`a $r$,
puis $c$ elle m\^eme de $r$ \`a $\sigma(p)$. On a alors
$$
\mathrm{longueur}(c)=\mathrm{longueur}(\tilde c)\ge\frac{\pi-\delta}{\sqrt 2}$$
d'apr\`es la Proposition  \ref{minsing}.
 \end{proof}

\section{Fin de la preuve}
\label{fin}

\subsection{Calcul des systoles dans le cas plat et dans le cas singulier}
\begin{lem}\label{typeb1}
Pour toute m\'etrique plate $g$ sur une vari\'et\'e $N$ de type $B_1$, on a 
$$\frac{(\mathrm{sys}(N,g))^3}{\mathrm{vol}(N,g)}\le\frac{2}{\sqrt{3}},$$
avec \'egalit\'e si et seulement si le r\'eseau associ\'e 
admet une base $(a_1,a_2,a_3)$ telle que $\vert{a_1}\vert=2\vert{a_2}\vert=2\vert{a_3}\vert$,
l'angle des vecteurs $a_1$ et $a_2$ \'etant  \'egal \`a $\pi/3$. Il existe sur $N$ une m\'etrique $g_0$, continue sur $N$ et lisse en 
dehors d'une hypersurface, localement isom\'etrique \`a $S^2\times\mathbb{R}$ l\`a o\`u elle est lisse, de quotient systolique strictement sup\'erieur \`a  $2/\sqrt{3}$.
\end{lem}
\begin{proof}[Preuve]
Dans le cas plat, le volume est \'egal \`a $\frac{1}{2}\det(a_1,a_2)\vert a_3\vert$, la systole
\`a $\inf\{\vert a_3\vert, s\}$, o\`u $s$ est la systole du tore plat de 
dimension $2$ d\'efini par le r\'eseau de base $\frac{a_1}{2},a_2$.
Normalisons de fa\c con que $|a_3|=1$. Alors le quotient systolique
vaut 
$$
\frac{2s^3}{\det(a_1,a_2)}\quad\hbox{si $s\le 1$, et}\quad \frac{2}{\det(a_1,a_2)}
\quad\hbox{si $s\ge 1$,}
$$
ce qui donne la premi\`ere partie d'apr\`es les r\'esultats classiques sur les r\'eseaux en
dimension $2$.

Pour la construction de $g_0$, on part de la m\'etrique sur la bouteille de Klein singuli\`ere
vue en \ref{sing},
dont la systole  est \'egale \`a $\pi$ et le volume \`a $2\pi\sqrt 2$,
et on proc\`ede par suspension :
d'apr\`es \ref{plat3} et  \ref{kleinbavard}  
le quotient de $(\mathbf{K},b)\times\RR$ par 
$<(r_\alpha,t_d)>$ 
est une vari\'et\'e  hom\'eomorphe \`a une vari\'et\'e de type $B_1$. 
Son volume est  \'egal \`a $2\pi\sqrt 2 d$, et sa systole \`a
$$
\inf\big\{\pi, \sqrt{ \mathbf{d}(r_\alpha)^2+d^2}\big\}.$$
D'apr\`es le Corollaire \ref{dklein}, $\mathbf{d}(r_\alpha)$ est \'egal \`a
$$ \frac{\alpha}{\sqrt{2}}\quad\hbox{si $\alpha \le \pi(2-\sqrt{2})$, et}\quad \pi-\alpha \quad\hbox{si $\pi(2-\sqrt{2}) < \alpha < \pi$} $$
Pour obtenir le meilleur quotient systolique en faisant cette construction, c'est-\`a-dire en faisant varier $\alpha$ et $d$,
 on fixe $\mathrm{sys}=\pi$, et on essaie de minimiser $d$. Pour le faire il faut maximiser $\mathbf{d}(r_\alpha)$, donc choisir
$\alpha=\pi(2-\sqrt 2)$ et  $d=\pi(2\sqrt 2 -2)^\frac{1}{2}$.

Pour ces valeurs de $\alpha$ et $d$ on a $\frac{\mathrm{sys}^3}{\mathrm{Vol}}=\frac{\pi}{4(\sqrt{2}-1)^{1/2}}>\frac{2}{\sqrt{3}}$
\end{proof}

\begin{lem}\label{typeb34}
Pour toute m\'etrique plate $(N,g)$ de type $B_3$ ou
$B_4$, on a 
$$\frac{\mathrm{sys}(N,g)^3}{\mathrm{vol}(N,g)}\le1,$$ avec \'egalit\'e si et seulement si 
 $\vert{a_1}\vert=\vert{a_2}\vert=2\vert{a_3}\vert$. Dans les deux cas, il existe sur $N$ une m\'etrique $g_0$
 continue,  lisse et localement isom\'etrique \`a $S^2\times\RR$ sur le com\-pl\'e\-mentaire d'une hypersurface, de quotient systolique strictement sup\'erieur \`a $1$.
\end{lem}
\begin{proof}[Preuve]
Dans les deux cas, le r\'eseau associ\'e est orthogonal et d'indice $4$ dans $\Gamma$,  si bien 
que 
$$\mathrm{vol}(N,g)=|a_1||a_2||a_3|/4\quad\hbox{et}\quad \mathrm{sys}(N,g)=\inf\{|a_1|/2, |a_2|/2,|a_3|\},$$ 
d'o\`u la premi\`ere partie.
Pour la deuxi\`eme partie, on proc\`ede par suspension, en partant de la bouteille
$(\mathrm{K},b)$ comme pr\'ec\'edemment, et en quotientant
$K\times\RR$ par le groupe $<(S_1,t_\pi)>$ (pour le type $B_3$)
ou $<(S_2,t_\pi)>$  (type $B_4$). Comme $S_1$ et $S_2$ ont des points fixes, la systole reste \'egale \`a $\pi$. Le volume vaut $2\pi^2\sqrt 2$, et le quotient
$(\mathrm{sys})^3/\mathrm{vol}$ est \'egal \`a $\frac{\pi}{2\sqrt 2}$.
\end{proof}

\begin{lem}\label{typeb2}
Pour toute m\'etrique plate $g$ sur une vari\'et\'e riemannienne $N$ de type $B_2$ on a 
$$\frac{\mathrm{sys}(N,g)^3}{\mathrm{vol}(N,g)}\le\frac{8}{\sqrt{39}}$$ 
L'\'egalit\'e est r\'ealis\'ee si et seulement si

i) $\vert{a_1}\vert=\vert{a_2}\vert=8d$, o\`u $d$ est
 la distance entre les plans $P_1$ et $P_2$ des sym\'etries gliss\'ees qui engendrent $\Gamma$ (cf.~\ref{plat3}) ;
 
 ii)
l'angle $\alpha$  des vecteurs $a_1$ et $a_2$ est \'egal \`a $\arccos(-15/24)$.

\end{lem}
\begin{proof}[Preuve]
On a $(N,g)=\mathbb{R}^3/\Gamma$, o\`u $\Gamma$ est d\'ecrit dans la d\'efinition de $B_2$. On note $\Lambda$ le r\'eseau associ\'e \`a $\Gamma$ et $L$ le r\'eseau $2$-dimensionnel engendr\'e par $(a_1,a_2)$. Le volume est \'egal \`a $\det(a_1,a_2){d}$. Dans la suite on identifie un r\'eseau $\Lambda$ avec ses \'elements. On sait que
$$\mathrm{sys}(N,g)=\inf_{\gamma \in \Gamma}\mathbf{d}(\gamma).$$
En regroupant les \'el\'ements de $\Gamma$ en trois cat\'egories, les sym\'etries gliss\'ees de vecteur $a_1+v$ ou $a_2+v$ ($v\in L$),
 et les \'el\'ements de $\Lambda$ c'est \`a dire les translations pures de $\Gamma$, on obtient
$$\mathrm{sys}(N,g)=\inf\{|a_1/2+\Lambda|,|a_2/2+\Lambda|,|\Lambda|\}.$$
En faisant intervenir le r\'eseau $L$, on obtient
$$\mathrm{sys}(N,g)=\inf\{|a_1/2+L|,|a_2/2+L |,  4d,\sqrt{|\frac{(a_1+a_2)}{2}+L|^2+4d^2},|L|\}$$

a) si $(a_1,a_2)$ est une base minimale du r\'eseau $L$ (les vecteurs $a_1$ et $a_2$ sont non colin\'eaires et minimaux dans $L$), alors 
$$\mathrm{sys}(N,g)=\inf\{|a_1|/2,|a_2|/2,4d\},\quad\mathrm{vol}(N,g)=|a_1||a_2|d\sin\alpha ,$$
d'o\`u 
$\frac{\mathrm{sys}^3}{\mathrm{vol}}\le 2/\sqrt{3}$ avec \'egalit\'e si $|a_1|=|a_2|=8d$ et $\alpha=\pi/3$

b) si $(a_1,a_1+a_2)$ est une base minimale de $L$ et $|a_1|\le|a_1+a_2|$ alors 
$$\mathrm{sys}(N,g)=\inf\{|a_1|/2,4d\}\ \hbox{et}\  \mathrm{vol}(N,g)=|a_1||a_2|d\sin\alpha$$ 
et on se ram\`ene au cas pr\'ec\'edent.

c) si $(a_1,a_1+a_2)$ est une base minimale de $L$ et $|a_1|\ge|a_1+a_2|$ alors  
$$\mathrm{sys}(N,g)=\inf\{|a_1|/2,4d,\sqrt{\frac{(a_1+a_2)^2}{4}+4d^2},|a_1+a_2|\}.$$
On adopte les normalisations d'usage concernant les r\'eseaux de dimension $2$, en rapportant
le plan \`a une base orthonorm\'ee $(u,v)$ telle que
$u=a_1+a_2$, $a_1=xu+yv$. Ainsi, $x^2+y^2\ge1$,  $\vert x\vert\le1/2$,
et $$\mathrm{vol}(N,g)=\frac{1}{2}\det(a_1,a_2,\frac{a_1+a_2}{2}+2dk)=dy.$$
On discute alors suivant les valeurs de $d$.

\begin{itemize}
\item
si $0<d\le1/8$ alors la systole vaut $4d$ et 
$$\frac{\mathrm{sys}^3}{\mathrm{vol}}=\frac{64d^2}{y}\le \frac{1}{y}\le\frac{2}{\sqrt{3}}$$
\item si $1/8\le d\le \frac{1}{4\sqrt{3}}$ alors $\mathrm{sys}=\inf\{\frac{\sqrt{x^2+y^2}}{2},4d\}$ et deux cas se pr\'esentent :

\begin{itemize}
\item
Premier cas  : $\frac{\sqrt{x^2+y^2}}{2} \ge 4d$ . Alors $\mathrm{sys}=4d$ et $y\ge \frac{\sqrt{256d^2-1}}{2}$, d'o\`u
$$\frac{\mathrm{sys}^3}{\mathrm{vol}}=\frac{64d^2}{y}\le \frac{128d^2}{\sqrt{256d^2-1}}\le \frac{8}{\sqrt{39}}$$ 
L'\'egalit\'e est r\'ealis\'ee si $d=\frac{1}{4\sqrt{3}}$, $x=1/2$ et $y=\frac{\sqrt{13}}{2\sqrt{3}}$
\item
Second cas : $\frac{\sqrt{x^2+y^2}}{2} \le 4d$. Alors $\mathrm{sys}=\frac{\sqrt{x^2+y^2}}{2}$ et 
le quotient systolique est \'egal \`a $\frac{(x^2+y^2)^{3/2}}{8yd}$. Il faut maximiser cette fonction sur  le domaine
$$\{(x,y,d) :x^2+y^2\ge1, \ x\le 1/2,\ 1/8\le d\le 1/4\sqrt{3}, \ 4d \ge \frac{\sqrt{x^2+y^2}}{2} \}$$
Une \'etude directe donne le m\^eme maximum qu'auparavant (m\^emes valeurs de $d$, $x$, $y$).
\end{itemize}
\item si
$\frac{1}{4\sqrt{3}} \le d\le \frac{\sqrt{3}}{4}$ alors $\mathrm{sys}(N,g)=\inf\{\frac{\sqrt{x^2+y^2}}{2},\frac{\sqrt{16d^2+1}}{2}\}$


\begin{itemize}
\item Premier cas: $\frac{\sqrt{x^2+y^2}}{2} \ge \frac{\sqrt{16d^2+1}}{2}$. Alors $$\frac{\mathrm{sys}^3}{\mathrm{vol}}=\frac{(16d^2+1)^{3/2}}{8yd}.$$ 
Mais  $y^2 \ge 16d^2+3/4$ d'o\`u $\frac{\mathrm{sys}^3}{\mathrm{vol}}\le 8d \frac{(16d^2+1)^{3/2}}{\sqrt{16d^2+3/4}}$.
Sur l'intervalle consid\'er\'e, cette fonction 
atteint son maximum en $d=\frac{1}{4\sqrt{3}}$. 
On retombe sur
$8/\sqrt{39}$.

\item
Second cas: $\frac{\sqrt{x^2+y^2}}{2} \le \frac{\sqrt{16d^2+1}}{2}$.
Alors  $$\frac{\mathrm{sys}^3}{\mathrm{vol}}
=\frac{(x^2+y^2)^{3/2}}{8yd}.$$
Comme pr\'ec\'edemment, le maximum est r\'ealis\'e en un point de la fronti\`ere du do\-maine de d\'e\-fi\-ni\-tion. Il est atteint lorsque $d=1/4\sqrt{3}$, $x=1/2$ et $4d=\frac{\sqrt{x^2+y^2}}{2}$ : c'est toujours $\frac{8}{\sqrt{39}}$.
\end{itemize}

\item si $d\ge \sqrt{3}/4$ alors $\mathrm{sys}=\inf \{1,\frac{\sqrt{x^2+y^2}}{2}\}$
\begin{itemize}
\item Premier cas: si $1 \le \frac{\sqrt{x^2+y^2}}{2}$ alors $\frac{\mathrm{sys}^3}{\mathrm{vol}}=\frac{1}{yd} \le \frac{2}{d\sqrt{15}}\le\frac{8}{\sqrt{45}}$
\item
Second cas: si $1 \ge \frac{\sqrt{x^2+y^2}}{2}$ 
alors  $\frac{\mathrm{sys}^3}{\mathrm{vol}}  =\frac{\sqrt{x^2+y^2}}{8yd} \le \frac{(x^2+y^2)^{3/2}}{2\sqrt{3}y}$.
Dans le domaine $\{4 \ge x^2+y^2\}$, cette fonction  atteint son maximum lorsque $x^2+y^2=4$. Ce maximum est aussi \'egal \`a $\frac{8}{\sqrt{45}}$.
\end{itemize}
\end{itemize}
\end{proof}

\begin{lem}\label{typeb2sing}
Si $N$ est une vari\'et\'e de type $B_2$, il existe sur $N$ une m\'etrique $g_0$ continue et lisse en dehors d'une hypersurface, localement iso\-m\'e\-trique \`a $S^2 \times \RR$, dont le quotient systolique est strictement sup\'erieur \`a $8/\sqrt{39}$.
\end{lem}

\begin{proof}[Preuve]
La construction de $g_0$ se fait comme pr\'ec\'edemment en partant de la bouteille de Klein singuli\`ere de Bavard. On fait le quotient de  $(\mathbf{K},b) \times \RR $ par $(T_\delta,t_d)$, o\`u $(T_\delta,t_d)$ est la transformation $T_\delta$ d\'ej\`a d\'efinie sur la composante $(\mathbf{K},b)$ et la translation de vecteur $d$ sur $\RR$. Cela  donne bien une vari\'et\'e hom\'eomorphe \`a une vari\'et\'e de type $B_2$, munie d'une m\'etrique singuli\`ere isom\'etrique localement \`a $S^2 \times \RR$ l\`a o\`u elle est lisse. On normalise en supposant que la systole vaut $ \pi$.
Cela impose l'in\'egalit\'e  $\big(\mathbf{d}(T_\delta,t_d)\big)^n \geq \pi$ pour tout $n$. Le volume est \'egal \`a $2\pi \sqrt{2}d$. Pour obtenir le plus grand quotient systolique il faut que $d$ soit minimal. On commence par \'etudier l'in\'equation $\mathbf{d}(T_\delta,t_d) \geq \pi$, qui est \'equivalente \`a 
$(\mathrm{d}(T_\delta))^2+d^2 \geq \pi^2$. La valeur maximale de $\mathrm{d}(T_\delta)$ est obtenue lorsque $\delta$ v\'erifie $\frac{\cos \delta - 1}{2}=\cos \frac{\pi-\delta}{\sqrt{2}}$. On note $\delta_0$ la solution de cette \'equation.
 Le $d$ minimal est obtenu lorsque l'\'egalit\'e est r\'ealis\'ee dans l'in\'equation pr\'ec\'edente, on le note $d_0$.  
on a bien 
$\mathbf{d}\big((T_{\delta_0},t_{d_0})^n\big)\geq \pi$ pour tout $n$ : comme $d_0$ est de l'ordre de $2,641$, $nd$ est sup\'erieur \`a $\pi$ pour $n \geq 2$. Ainsi, ces deux constantes donnent le plus grand quotient systolique parmi ceux de toutes les constructions similaires. Il est \'egal \`a $\frac{\pi^2}{2\sqrt{2} d_0} >8/ \sqrt{39}$. 
\end{proof}

Le tableau suivant permet de faire la comparaison entre les plus grands  quotients systoliques des m\'etriques  plates ($\tau$(plates)) et les 
les plus grands quotients systoliques des  m\'etriques singuli\`eres que nous venons de construire ($\tau$(singuli\`eres))  sur chaque type de vari\'et\'e de Bieberbach non-orientables de  dimension $3$. Ici, $d_0$ d\'esigne la valeur de $d$ obtenue \`a la fin de la preuve
ci-dessus.

\begin{table}[htbp]\renewcommand{\arraystretch}{1.2}
\begin{center}
\begin{tabular}[b]{|*{5}{c|}}
   \hline
type    &$\tau$(plates)&valeur approch\'ee&$\tau$(singuli\`eres)&valeur  approch\'ee\\
   \hline
   $B_1$&$\frac{2}{\sqrt{3}}$&$\thickapprox  1,154$&$\frac{\pi}{4(\sqrt{2}-1)^{\frac{1}{2}}}$&$\thickapprox 1,220$\\
   \hline
   $B_2$&$\frac{8}{\sqrt{39}}$&$\thickapprox  1,281$&$\frac{\pi2}{2d_0\sqrt{2}}$&$\thickapprox 1,321$\\
   \hline
   $B_3$&$1$&$1$&$\frac{\pi}{2\sqrt{2}}$&$\thickapprox 1,110$\\
   \hline
   $B_4$&$1$&$1$&$\frac{\pi}{2\sqrt{2}}$&$\thickapprox 1,110$\\
   \hline
   \end{tabular}\end{center}
\end{table}

\subsection{Des m\'etriques singuli\`eres aux m\'etriques lisses}
\label{continuite}
Les m\'etriques que nous venons de construire sont seulement continues.
Pour obtenir des m\'etriques lisses ayant les m\^emes propri\'et\'es, on utilise des arguments de continuit\'e
et d'approximation.

Si $M$ est une vari\'et\'e compacte, l'espace des sections continues de n'importe quel fibr\'e tensoriel sur $M$
peut \^etre muni de la topologie \emph{compacte ouverte}, cf. \cite{hirsch}, ch. 2. 
Pour cette topologie,
le sous-espace $\mathcal{M}$ des m\'etriques riemanniennes est ouvert dans l'espace $C^0(S^2M)$
des $2$-tenseurs sym\'etriques. Pour la topologie induite sur  $\mathcal{M}$, on peut v\'erifier que
tout ouvert est une r\'eunion d'ensembles de la forme
$$
A(g,\epsilon)=\big\{h\in\mathcal{M}, \forall x\in M,\, (1-\epsilon)g_x<h_x<(1+\epsilon)g_x\big\},$$
o\`u $g\in\mathcal{M}$ et $\epsilon\in]0,1[$, ces in\'egalit\'es \'etant prises au sens des formes quadratiques sur $T_xM$.
 
La propri\'et\'e suivante est folklorique. Faute de pouvoir en donner une r\'ef\'erence \'ecrite, nous
la d\'emontrons.

\begin{theo}
L'application $h\mapsto \mathrm{sys}(h)$ est continue sur $\mathcal{M}$ muni de la topologie
compacte ouverte.
\end{theo}
\begin{proof}[Preuve]
Il suffit de travailler avec un $A(g,\epsilon)$.
Pour $h\in A(g,\epsilon)$, la systole peut \^etre r\'ealis\'ee par une courbe
param\'etr\'ee $c :[0,1]\rightarrow M$
prise dans l'ensemble 
$$\mathcal{C}=\cup_{h\in A(g,\epsilon)}\mathcal{C}_h,$$
o\`u $\mathcal{C}_h$ est l'ensemble des lacets $c:[0,1]\rightarrow M$
non contractiles, $C^1$ par morceaux, param\'etr\'es propor\-tionnellement \`a la longueur d'arc pour $h$,
et de $h$-longueur in\-f\'erieure ou \'egale \`a 
$2\mathrm{diam}(M,h)$.

Dans ces conditions, l'ensemble des longueurs des courbes de $\mathcal{C}$, prises pour toutes les m\'etriques de
 $A(g,\epsilon)$, est major\'e par $2(1+\epsilon)\mathrm{diam}(M,g)$, et l'ensemble
$$\bigcup\big\{\big(c(t),\dot c^{\pm}(t)\big)\big\}\quad\hbox{avec $c\in\mathcal{C}$, $t\in[0,1]$}$$
	est relativement compact dans $TM$. Il en r\'esulte que la famille
(index\'ee par $\mathcal{C}$) de fonctions qui \`a une m\'etrique $h$
associe la longueur de $c$ pour $h$ est \emph{\'equicontinue} sur $A(g,\epsilon).$
Mais la borne inf\'erieure d'une famille de fonctions \'equicontinue est continue.
\end{proof}

Enfin (voir encore \cite{hirsch}), l'espace des sections lisses de n'importe quel fibr\'e vectoriel de base compacte est
dense pour la topologie compacte ouverte dans l'espace des sections continues. Dans la situation qui nous int\'eresse,
il suffit en fait d'approcher uniform\'ement la fonction $f$ par une fonction lisse.

\section{Remarques diverses}

Une m\'etrique de la forme
$d\phi^2 +f^2(\phi)d\theta^2$ est localement conforme \`a une m\'etrique plate, m\^eme si $f$ est seulement continue~: il suffit de mettre
$f^2$ en facteur. Le m\^eme proc\'ed\'e montre
de plus que la m\'etrique de C.~Bavard est globalement conforme \`a une m\'etrique plate sur la
bouteille de Klein. Cette m\'etrique est facile \`a expliciter (voir \cite{bavard2} pour les d\'etails).

Ici, comme nous travaillons en dimension $3$, la situation est diff\'erente. Le mod\`ele local est de la forme 
$d\phi^2 +f^2(\phi)d\theta^2+dt^2$.
Il s'av\`ere que si $g$ est une m\'etrique riemannienne $C^3$ en dimension $2$,
$g+dt^2$ est conform\'ement plate si et seulement si $g$ est \`a courbure constante.
Par cons\'equent, les m\'etriques que nous avons introduites en $\ref{fin}$ sont bien conform\'ement plates
\emph{en dehors de la  surface singuli\`ere}. Il ne faut pas s'attendre \`a ce qu'elles soient conform\'ement plates
au voisinage d'un point d'une surface singuli\`ere, mais la m\'etrique \'etant seulement continue, un argument sp\'ecifique est
n\'ecessaire.

\begin{prop}
La m\'etrique $d\phi^2 +f^2(\phi)d\theta^2+ dt^2$,
o\`u $f$ est la fonction donn\'ee dans (\ref{sph}),
n'est pas conform\'ement plate.
\end{prop}

\begin{proof}[Preuve]
Soit $p$ un point de l'hypersurface singuli\`ere $S$ d'\'equation $\phi=
\phi_0$, et $h:U\rightarrow\RR^3$ un diff\'eomorphisme conforme sur un ouvert de $\RR^3$.
Quitte \`a diminuer $U$, on peut le supposer invariant par la transformation
$s: (\phi,\theta, t)\mapsto (2\phi_o-\phi,\theta,t)$,
qui est une isom\'etrie de $U$ fixant $S$.
D'apr\`es le th\'eor\`eme de Liouville (voir par exemple \cite{kp}, p.12),
la conjugu\'ee de $s$ par $h$ est la restriction d'une transformation de M\"obius.
Mais une telle transformation ne peut exister : l'ensemble de ses points fixes
contient $h(S\cap U)$, qui devrait \^etre une portion de $2$-sph\`ere. 
Pour voir que c'est impossible, on introduit l'application conforme
$h_1$ de $U^\prime =U\cap\{\phi<\phi_o\}$
dans $\RR^3$ donn\'ee par
$$h_1(\phi,\theta, t)=(e^t\cos\phi\cos\theta,e^t\cos\phi\sin\theta,e^t\sin\phi)$$
(comparer \`a \cite{kp}, p.76--77).
Toujours d'apr\`es le th\'eor\`eme de Liouville, il existe une transformation de M\"obius $h_2$ 
telle qu'en restriction \`a $U^\prime$
on ait $h=h_2\circ h_1$.
Par continuit\'e, cette propri\'et\'e reste vraie sur $S\cap U$.
Ainsi, l'ensemble des points fixes de la transformation de M\"obius $hsh^{-1}$ contient $h_2\big(h_1(S\cap U)\big)$. On aboutit \`a  la contradiction cherch\'ee, car
$h_1(S)$ est une portion de c\^one de r\'evolution.
\end{proof}

Pour plus de d\'etails sur les r\'esultats de g\'eom\'etrie conforme utilis\'es ici, voir
les contributions de R. Kulkarni et du second auteur dans \cite{kp}.

Notons enfin
que pour celles de ces m\'etriques qui donnent le meilleur quotient systolique, la vari\'et\'e est encore recouverte
par les g\'eod\'esiques systoliques. Il n'y a donc pas d'obstacle ``\'evident(!)'' \`a ce que ces m\'etriques
soient optimales, m\^eme si nous sommes loin de pouvoir le montrer. 

\medskip
\noindent
\textbf{Remerciements.} C'est un plaisir de remercier le referee pour ses critiques constructives.


\begin{thebibliography}{99}



\bibitem{bavard}  Bavard, C., In\'egalit\'e isosystolique pour la bouteille de Klein, Math. Ann. 274, 439--441(1986)

\bibitem{bavard2} Bavard, C., In\'egalit\'es isosystoliques conformes  pour la bouteille de Klein, Geom. Dedicata 27, 349--355 (1988), 


\bibitem{berger} Berger, M., Quelques probl\`emes de g\'eom\'etrie riemannienne ou deux variations sur les espaces
compacts sym\'etriques de rang $1$, L'Ens.Math. (2) 16  73--96 (1970).

\bibitem{bki} Berger, M., Systoles et applications selon Gromov, S\'eminaire N. Bourbaki, 
expos\'e 771, Ast\'erisque 216, 279--310 (1993).

\bibitem{bbi} Burago D., Burago, Y.D., Ivanov, S., A course in metric geometry,
Graduate studies in Mathematics (33), Amer. Math. Soc., Providence, R.I. 2001.

\bibitem{charlap} Charlap, L.S., Bieberbach Groups and Flat Manifolds, Universitext, Berlin 1986.



\bibitem{ghl} Gallot, S., Hulin, D., Lafontaine, J., Riemannian Geometry, 3rd edition, Springer,
Berlin Heidelberg 2004.

\bibitem{gromov} Gromov, M., Filling Riemannian manifolds, J. Diff. Geom. 18, 1--147(1983)

\bibitem{gromov1} Gromov M., Systoles and intersystolic inequalities,  in : Besse, A.L. (ed.),
Actes de la table ronde de g\'eom\'etrie diff\'erentielle en l'honneur de 
Marcel Berger, Soci\'et\'e Math\'ematique de France, S\'eminaires et Congr\`es no. 1, p. 291--362.

\bibitem{hirsch}  Hirsch  M., Differential Topology, Springer, Berlin Heidelberg 1976

\bibitem{katz} Katz, M.G, Systolic Geometry and Topology, Math. Surveys and Monographs137,
Amer. Math. Soc., Providence, R.I. 2007.

\bibitem{kp} Kulkarni, R.S., Pinkall, U. (Eds), Conformal Geometry, Aspects of Mathematics vol.12, Vieweg,
Braunschweig 1988 


\bibitem{pu}   Pu, P.M.,  Some inequalities in certain non-orientable riemannian manifolds. Pacific J.Math.2, 55--71(1952)

\bibitem{sakai} Sakai, T., A proof of the isosystolic inequality for the Klein bottle, Proc. Amer. Math. Soc. 104, 589--590 (1988)

\bibitem{th} Thurston, W.P., Three-Dimensional Geometry and Topology, edited by S. Levy, Princeton University Press,
Princeton (1997)

\bibitem{wolf} Wolf, J.A., Spaces of constant curvature, Publish or Perish, Boston 1974

\end{thebibliography}
\end{document}